\documentclass[11pt]{article}
\usepackage[cm]{fullpage}

\usepackage{latexsym}
\usepackage{amsfonts}
\usepackage{amsmath}
\usepackage{amssymb}
\usepackage{amscd}

\newtheorem{theorem}{Theorem}[section]

\newtheorem{remark}[theorem]{Remark}
\newtheorem{corollary}[theorem]{Corollary}

\begin{document}

\title{\normalsize{ON POTENTIAL FUNCTION OF GRADIENT STEADY RICCI SOLITONS}}
\author{\scriptsize{PENG WU}}
\date{}
\maketitle

\begin{abstract}
In this paper we study potential function of gradient steady Ricci
solitons. We prove that infimum of potential function decays
linearly; in particular, potential function of rectifiable gradient
steady Ricci solitons decays linearly. As a consequence, we show
that a gradient steady Ricci soliton with bounded potential function
must be trivial, and no gradient steady Ricci soliton admits
uniformly positive scalar curvature.
\end{abstract}

\section{Introduction}
A complete Riemannian metric $g$ on a smooth manifold $M^n$ is
called a \textit{gradient Ricci soliton} if there is a smooth
function $f:M^n\rightarrow \mathbb{R}$ such that
\begin{equation}
\begin{split}
Ric+\text{Hess} f &= \lambda g,
\end{split}
\end{equation}
for some constant $\lambda$. The function $f$ is called a
\textit{potential function} for $g$. When $\lambda>0$ the Ricci
soliton is \textit{shrinking}, when $\lambda=0$ it is
\textit{steady}, and when $\lambda<0$ it is \textit{expanding}. When
$f$ is constant the gradient Ricci soliton is simply an Einstein
manifold. Thus Ricci solitons are natural extensions of Einstein
metrics, an Einstein manifold with constant potential function is
called a trivial gradient Ricci soliton. Gradient Ricci solitons
play an important role in Hamilton's Ricci flow as they correspond
to self-similar solutions, and often arise as singularity models.
They are also related to smooth metric measure spaces, since
equation (1) is equivalent to $\infty$-Bakry-Emery Ricci tensor
$Ric_f=0$. In physics, a smooth metric space $(M^n,\ g,\
e^{-f}dvol)$ with $Ric_f=\lambda g$ is called quasi-Einstein
manifold. Therefore it is important to study geometry and topology
of gradient Ricci solitons and their classifications. Recently there
have been a lot of work on gradient solitons, see
\cite{CCZ,Chu,Der,ENM,MS1,Nab,NW1,NW2,PW1,PW2,PW3,PRS} for example;
and \cite{Cao1,Cao2} for excellent surveys.


The growth of the potential function of gradient Ricci solitons has
been an interesting problem. For gradient shrinking solitons, H.-D.
Cao and D. Zhou \cite{CZ} proved that the potential function grows
quadratically, and based on this estimate they proved that gradient
shrinking Ricci solitons have at most Euclidean volume growth.


In \cite{H} R. Hamilton proved the following identity for gradient
Ricci solitons
\begin{equation}
\begin{split}
R+|\nabla f|^2-2\lambda f=\Lambda.
\end{split}
\end{equation}
where $R$ is the scalar curvature of $(M^n,g)$, and $\Lambda$ is a
constant. For gradient steady Ricci solitons, equation (2) becomes
\begin{equation*}
R+|\nabla f|^2 = \Lambda.
\end{equation*}
where $\Lambda$ is a constant. B.-L. Chen \cite{Chen} proved that
the scalar curvature $R\geq 0$, hence $|\nabla
f|\leq\sqrt{\Lambda}$, and the potential function has at most linear
growth. It is natural to ask that whether potential function of a
gradient steady Ricci soliton grows linearly at infinity.

Recently H.-D. Cao and Q. Chen \cite{CC} partially confirmed this
under additional assumption of positive Ricci curvature, they proved
that if a gradient steady Ricci soliton has positive Ricci curvature
and scalar curvature attains its maximum at some point, then the
potential function decays linearly.

In general one cannot expect potential function to grow or decay
linearly along all directions at infinity, because of the product
property: the product of any two gradient steady Ricci solitons is
also a gradient steady Ricci soliton. Consider for example
$(\mathbb{R}^2,g_0,f)$, where $g_0$ is the standard Euclidean
metric, $f(x_1,x_2)=x_1$. $f$ is constant along $x_2$ direction, so
without additional conditions, $f$ may not have linear growth at
infinity. We prove that though the potential function may not be
linear, its infimum must decay linearly,

\begin{theorem}
Let $(M^n,g,f)$ be a gradient steady Ricci soliton with $R+|\nabla
f|^2=\Lambda$. For any $x\in M^n$, there exists $r_0>0$, such that
for any $r\geq r_0$,
\begin{equation*}
\begin{split}
-\sqrt{\Lambda}r \leq \inf_{y\in
\partial B_r(x)} f(y)-f(x) &\leq -\sqrt{\Lambda}r + \sqrt{2n\sqrt{\Lambda}}(\sqrt{r}+1).
\end{split}
\end{equation*}
Therefore, a gradient steady Ricci soliton with bounded potential
function must be trivial.
\end{theorem}

\begin{remark}
The author was told that Ovidiu Munteanu and Natasa Sesum \cite{MS2}
also got the same result using a different method.
\end{remark}
As a consequence, we have
\begin{corollary}
Let $(M^n,g,f)$ be a complete noncompact gradient steady Ricci
soliton. Then for any $x\in M^n$,
\begin{equation*}
\begin{split}
\limsup_{y\in B_r(x), r\rightarrow \infty} |\nabla f|(y) &=
\sqrt{\Lambda},\ \ \ \liminf_{y\in B_r(x), r\rightarrow \infty} R(y)
= 0.
\end{split}
\end{equation*}
In another word, no gradient steady Ricci soliton admits uniformly
positive scalar curvature.
\end{corollary}

In \cite{PW1} P. Petersen and W. Wylie introduced the concept of
rectifiable Ricci soliton and studied their rigidity. A gradient
Ricci soliton is called rectifiable if its potential function
$f=f(r)$, where $r$ is the distance function. In this case
$\inf_{y\in \partial B_r(x)}=f(r)$, by the same argument for proof
of Theorem 1.1 we show that its potential function decays linearly,

\begin{theorem}
Let $(M^n,g,f)$ be a complete noncompact rectifiable gradient steady
Ricci soliton. Then there exists $r_0>0$, such that when
$r=d(x,\cdot)\geq r_0$,
\begin{equation*}
\begin{split}
-\sqrt{\Lambda}r \leq  f(r)-f(0) &\leq -\sqrt{\Lambda}r +
\sqrt{n\sqrt{\Lambda}}(\sqrt{r}+1).
\end{split}
\end{equation*}
\end{theorem}

\begin{remark}
We in fact obtained an estimate of potential function for almost all
known gradient steady Ricci soliton examples, since almost all of
them are rectifiable, see for instance solitons constructed in
\cite{DW} and \cite{Ivey}.
\end{remark}

Our tool is $f$-volume comparison theorem for smooth metric measure
spaces that were developed by G. Wei and W. Wylie \cite{WW}. In
Section 2 we improve the $f$-volume comparison theorem for smooth
metric measure spaces with nonnegative Bakry-Emery Ricci tensor, and
apply to gradient steady Ricci solitons to obtain an upper bound of
the $f$-volume. In Section 3 we deduce a lower bound for the
$f$-volume and prove the main theorem.

\

\textbf{Acknowledgments.} The author would like to thank his
advisors Xianzhe Dai and Guofang Wei for their constant support and
help.

\section{Volume Comparison for Smooth Metric Measure Spaces}

In this section, we improve $f$-volume comparison theorem for smooth
metric spaces under the condition $f$ is at most linear. Recall that
a smooth metric measure space is triple $(M^n, g,
e^{-f}\text{dvol}_g)$, where $(M^n, g)$ is a Riemannian manifold,
$f$ is a smooth real valued function on $M$. smooth metric measure
spaces play an essential role in Perelman's work on the Ricci flow,
and they arise as smooth collapsed measured Gromov-Hausdorff limits.

The Ricci curvature of smooth metric measure space, which is called
Bakry-Emery Ricci curvature, is defined as
Ric$_f=\text{Ric}_g+\text{Hess}f$. The self-adjoint Laplacian with
respect to the weighted measure $e^{-f}\text{dvol}_g$ is $\Delta_f
=\Delta -\nabla_{\nabla f}$, the weighted (or $f$-)mean curvature is
defined as $m_f=m-\langle\nabla f,\nabla r\rangle=\Delta_f(r)$,
where $r$ is the distance function, and the weighted (or $f$-)volume
is defined as Vol$_f (B_r(x))=\int_{B_r(x)} e^{-f}\text{dvol}_g$.
Fix $x\in M^n$, under exponential polar coordinates around $x$ we
write the volume element
$\text{dvol}_g=\mathcal{A}(r,\theta)dr\wedge d\theta_{n-1}$, where
$d\theta_{n-1}$ is the standard volume element on the unit sphere
$S^{n-1}$. Let
$\mathcal{A}_f(r,\theta)=e^{-f}\mathcal{A}(r,\theta)$, it is easy to
check that $(\ln(\mathcal{A}_f(r,\theta)))'=m_f(r)$, and Vol$_f
(B_r(x))=\int_0^r\int_{S^{n-1}} \mathcal{A}_f(t,\theta)dtd\theta$.

Denote Vol$^n_H(B_r)$ be the volume of the radius $r$-ball in the
model space $M^n_H$. G. Wei and W. Wylie \cite{WW} proved the
following $f$-volume comparison theorem for smooth metric measure
spaces,

\begin{theorem} ($f$-Volume comparison) (Theorem 1.2 in \cite{WW}).\\
Suppose $(M^n,g,e^{-f}dvol)$ is a smooth metric measure space with
$Ric_f\geq (n-1)H$. Fix $x\in M$. If $|f|\leq \Lambda$. Then for
$R\geq r>0$ (assume $R\leq \pi\slash 4\sqrt{H}$ if $H>0$)
\begin{equation*}
\begin{split}
\frac{V_f(B_R(x))}{V_f(B_r(x))} &\leq
\frac{V_H^{n+4\Lambda}(B_R)}{V_H^{n+4\Lambda}(B_r)}.
\end{split}
\end{equation*}
\end{theorem}

\begin{remark} If $|\nabla f|\leq a$, without loss of generality,
assume $f(x)=0$, then $f(y)\leq aR$ for $y\in B_R(x)$, and Theorem
2.1 becomes
\begin{equation*}
\begin{split}
\frac{V_f(B_R(x))}{V_f(B_r(x))} &\leq
\frac{V_H^{n+4Ra}(B_R)}{V_H^{n+4Ra}(B_r)}.
\end{split}
\end{equation*}
\end{remark}

For our purpose, we concentrate on the case $H=0$. Denote
Vol$^n_{\mathbb{R}}(B_r)$ be the volume of the radius $r$-ball in
$\mathbb{R}^n$. When $f$ has at most linear growth, we obtain the
following:

\begin{theorem}
Suppose $(M^n,g,e^{-f}dvol)$ is a smooth metric measure space with
$Ric_f\geq 0$, $|\nabla f|\leq a$; and in addition
$\frac{f(y)-f(x)}{d(y,x)}\geq -a+\epsilon$ for $y\in
B_R(x)\backslash B_s(x)$ for $R\geq s>0$, then for $s<S<R, s<r<R$,
we have
\begin{equation*}
\begin{split}
\frac{V_f(B_R(x)\backslash B_r(x))}{V_f(B_S(x)\backslash B_s(x))}
&\leq \frac{V_{\mathbb{R}}^{n+R\bar a}(B_R\backslash
B_r)}{V_{\mathbb{R}}^{n+R\bar a}(B_S\backslash B_s)}.
\end{split}
\end{equation*}
where $\bar a=a -\epsilon^2/2a$. In particular, if
$\frac{f(y)-f(x)}{d(y,x)}\geq -a+\epsilon$ for $y\in B_R(x)$, then
for $r<R$, we have
\begin{equation*}
\begin{split}
\frac{V_f(B_R(x))}{V_f(B_r(x))} &\leq \frac{V_{\mathbb{R}}^{n+R\bar
a}(B_R)}{V_{\mathbb{R}}^{n+R\bar a}(B_r)}.
\end{split}
\end{equation*}
\end{theorem}

\begin{remark}
In Theorem 2.3 we improved the dimension of the model space from
$n+4ar$ to $n+\bar a r$, which is crucial in our proof of the main
theorem.
\end{remark}

To prove the $f$-volume comparison theorem, we first prove the
following $f$-mean curvature comparison theorem, and Theorem 2.3
follows directly from the argument in G. Wei and W. Wylie \cite{WW}
and Lemma 3.2 in \cite{Zhu}.

\begin{theorem}($f$-Mean Curvature Comparison).
Suppose $(M^n,g,e^{-f}dvol)$ is a smooth metric measure space with
$Ric_f\geq 0$. If $|\nabla f|\leq a$ along a minimal geodesic
segment from $x$, and $\frac{f(y)-f(x)}{d(y,x)}\geq -a+\epsilon$ for
$y\in\partial B_r(x)$, then
\begin{equation*}
\begin{split}
m_f(r) &\leq \bar a +\frac{n-1}{r}=m_{\mathbb{R}}^{n+\bar ar}(r).
\end{split}
\end{equation*}
along that minimal geodesic segment from $x$, where $\bar a = a
-\epsilon^2/2a$.
\end{theorem}
$\mathbf{Proof}$. From inequality (2.21) in \cite{WW}
\begin{equation*}
\begin{split}
m_f(r) &\leq \frac{n-1}{r} -\frac{2}{r}f(r) +\frac{2}{r^2}\int_0^r
f(t)dt.
\end{split}
\end{equation*}

Suppose $f(y)-f(x)=(-a+\epsilon)r$ for some $y\in\partial B_r(x)$,
we will maximize
\[-\frac{2}{r}f(y) +\frac{2}{r^2}\int_0^r
f(t)dt=\frac{2}{r^2}\int_0^r (f(t)-f(y))dt\] along a minimal
geodesic segment from $x$ to $y$.

Since $|\nabla f|\leq a$, along a minimal geodesic segment from $x$
to $y$, $f(t)-f(x)$ is bounded from above by
\begin{equation*}
F(t)=\left\{
\begin{aligned}
a t, & \ \ \ \ 0\leq t\leq \frac{\epsilon}{2a}r,\\
-a t +\epsilon r, & \ \ \ \ \frac{\epsilon}{2a}r \leq t\leq r,\\
\end{aligned}
\right.
\end{equation*}
Thus
\begin{equation*}
\begin{split}
\int_0^r (f(t)-f(y))dt &= \int_0^r [(f(t)-f(x))-(f(y)-f(x))]dt\\
&\leq \int_0^r (F(t)-(-a+\epsilon)t) dt\\
&= \frac{r^2}{2}(a-\frac{\epsilon^2}{2a}).
\end{split}
\end{equation*}

Hence if $\frac{f(y)-f(x)}{d(y,x)}\geq -a+\epsilon$ on $\partial
B_{r}(x)$, we obtain the following $f$-mean curvature comparison
\begin{equation*}
\begin{split}
m_f(r) &\leq \bar a + \frac{n-1}{r} = m_{\mathbb{R}}^{n+\bar a
r}(r).
\end{split}
\end{equation*}

\section{Proof of the Main Theorems}

Proof of Theorem 1.1. To prove that infimum of the potential
function decays linearly, we will derive an inequality on the left
hand side of the volume comparison in Theorem 2.3.

Taking trace of the soliton equation (1), we get $R+\Delta f=0$. Add
to equation (2) we obtain
\begin{equation}
-\Delta_f f=-\Delta f+|\nabla f|^2=\Lambda.
\end{equation}
Hence by maximal principle, $f$ has no local minimum.

Choose $x\in M,\ \delta>0$. Suppose $\inf \frac{f(y)-f(x)}{d(y,x)} =
-\sqrt{\Lambda}+\epsilon$ in $B_{r+\sqrt{r}+\delta}(x)\backslash
B_r(x)$, so $\bar
\Lambda=\sqrt{\Lambda}-\epsilon^2/2\sqrt{\Lambda}$. Since $f$ has no
local minimum, we have $\epsilon\leq\sqrt{\Lambda}$, and
$\bar\Lambda\geq \sqrt{\Lambda}/2$. Choose a smooth cut-off function
$\phi$, such that
\begin{equation*}
\phi(y)= \left\{
\begin{aligned}
1, & \hspace{1cm} y\in B_{r+\sqrt{r}}(x)\\
0, & \hspace{1cm} y\in M\backslash B_{r+\sqrt{r}+\delta}(x),
\end{aligned}
\right.
\end{equation*}
and $|\nabla \phi|\leq \frac{1+\delta}{\delta}$.

Integrate equation (3) in $B_{r+\sqrt{r}+\delta}(x)$ and apply
stokes formula,
\begin{equation*}
\begin{split}
\Lambda\int_{B_{r+\sqrt{r}+\delta}(x)}e^{-f}\phi \text{dvol}
&= -\int_{B_{r+\sqrt{r}+\delta}(x)} \Delta_f f\ e^{-f} \phi \text{dvol}\\
&= \int_{B_{r+\sqrt{r}+\delta}(x)\backslash B_{r+\sqrt{r}}(x)} \langle\nabla f, \nabla\phi\rangle e^{-f} \text{dvol}\\
& \leq \int_{B_{r+\sqrt{r}+\delta}(x)\backslash B_{r+\sqrt{r}}(x)}
|\nabla f||\nabla\phi| e^{-f} \text{dvol}\\
&\leq \frac{(1+\delta)\sqrt{\Lambda}}{\delta}
\int_{B_{r+\sqrt{r}+\delta}(x)\backslash B_{r+\sqrt{r}}(x)} e^{-f}
\text{dvol}.
\end{split}
\end{equation*}
Therefore
\begin{equation*}
\begin{split}
\Lambda\int_{B_{r+\sqrt{r}}(x)\backslash B_r(x)} e^{-f} \text{dvol}
&\leq \Lambda\int_{B_{r+\sqrt{r}+\delta}(x)} e^{-f}\phi \text{dvol}\\
&\leq \frac{(1+\delta)\sqrt{\Lambda}}{\delta}
\int_{B_{r+\sqrt{r}+\delta}(x)\backslash B_{r+\sqrt{r}}(x)} e^{-f}
\text{dvol}.
\end{split}
\end{equation*}
So we get
\begin{equation*}
\begin{split}
\frac{\delta\sqrt{\Lambda}}{1+\delta} &\leq
\frac{V_f(B_{r+\sqrt{r}+\delta}(x)\backslash
B_{r+\sqrt{r}}(x))}{V_f(B_{r+\sqrt{r}}(x)\backslash B_r(x))}.
\end{split}
\end{equation*}
On the other hand, by Theorem 2.3, we have
\begin{equation*}
\begin{split}
\frac{V_f(B_{r+\sqrt{r}+\delta}(x)\backslash
B_{r+\sqrt{r}}(x))}{V_f(B_{r+\sqrt{r}}(x)\backslash B_r(x))} &\leq
\frac{V_{\mathbb{R}}^{n+(r+\sqrt{r}+\delta)\bar\Lambda}(B_{r+\sqrt{r}+\delta})
-V_{\mathbb{R}}^{n+(r+\sqrt{r}+\delta)\bar\Lambda}(B_{r+\sqrt{r}})}
{V_{\mathbb{R}}^{n+(r+\sqrt{r}+\delta)\bar\Lambda}(B_{r+\sqrt{r}})
-V_{\mathbb{R}}^{n+(r+\sqrt{r}+\delta)\bar\Lambda}(B_{r})}\\
&= \frac{(r+\sqrt{r}+\delta)^{n+(r+\sqrt{r}+\delta)\bar\Lambda}
-(r+\sqrt{r})^{n+(r+\sqrt{r}+\delta)\bar\Lambda}}
{(r+\sqrt{r})^{n+(r+\sqrt{r}+\delta)\bar\Lambda}
-r^{n+(r+\sqrt{r}+\delta)\bar\Lambda}}.
\end{split}
\end{equation*}
Therefore,
\begin{equation*}
\begin{split}
\frac{\delta\sqrt{\Lambda}}{1+\delta} &\leq
\frac{(r+\sqrt{r}+\delta)^{n+(r+\sqrt{r}+\delta)\bar\Lambda} -
(r+\sqrt{r})^{n+(r+\sqrt{r}+\delta)\bar\Lambda}}
{(r+\sqrt{r})^{n+(r+\sqrt{r}+\delta)\bar\Lambda}
-r^{n+(r+\sqrt{r}+\delta)\bar\Lambda}}.
\end{split}
\end{equation*}
Divide both sides by $\delta$ and let $\delta\rightarrow 0$, we get
\begin{equation*}
\begin{split}
\sqrt{\Lambda} &\leq
\frac{\frac{n}{r+\sqrt{r}}+\bar\Lambda}{1-(1+\frac{1}{\sqrt{r}})^{-n-(r+\sqrt{r})\bar\Lambda}}.
\end{split}
\end{equation*}
Therefore,
\begin{equation*}
\begin{split}
\epsilon \leq
\sqrt{\frac{2n\sqrt{\Lambda}}{r+\sqrt{r}}+2\Lambda(1+\frac{1}{\sqrt{r}})^{-n-\sqrt{\Lambda}(r+\sqrt{r})/2}}.
\end{split}
\end{equation*}
Thus when $r$ is sufficiently large, $\epsilon\leq
\sqrt{\frac{2n\sqrt{\Lambda}}{r}}$, hence
\begin{equation*}
\begin{split}
\inf_{y\in B_{r+\sqrt{r}}(x)\backslash B_r(x)}
\frac{f(y)-f(x)}{d(x,y)} \leq
-\sqrt{\Lambda}+\sqrt{\frac{2n\sqrt{\Lambda}}{r}}.
\end{split}
\end{equation*}

Now suppose $\inf\frac{f(y)-f(x)}{d(y,x)}$ is attained by $y_0\in
B_{r_0}(x)$, for some $r\leq r_0\leq r+\sqrt{r}$. Let $z=\partial
B_r(x)\cap\gamma(t)$, where $\gamma(t)$ is the minimal geodesic from
$x$ to $y_0$, then
\begin{equation*}
\begin{split}
f(y_0)-f(x) &\leq  -\sqrt{\Lambda}r_0+\sqrt{\frac{2n\sqrt{\Lambda}}{r}}r_0,\\
f(z)-f(y_0) &\leq \sqrt{\Lambda}(r_0-r).
\end{split}
\end{equation*}
Therefore,
\begin{equation*}
\begin{split}
\inf_{y\in \partial B_r(x)} f(y)-f(x) &\leq f(z)-f(x)\leq
-\sqrt{\Lambda}r +\sqrt{2n\sqrt{\Lambda}}(\sqrt{r}+1).
\end{split}
\end{equation*}

\

\

Proof of Theorem 1.4. For rectifiable gradient steady Ricci
solitons, since $f(r)$ has no local minimum, it attains its maximum
at $r=0$, and is monotonically nonincreasing, hence the function
$F(t)$ in the proof of Theorem 2.5 can be replaced by
\begin{equation*}
F(t)=\left\{
\begin{aligned}
0, & \ \ \ \ 0\leq t\leq \frac{\epsilon}{a}r,\\
-a t +\epsilon r, & \ \ \ \ \frac{\epsilon}{a}r \leq t\leq r,\\
\end{aligned}
\right.
\end{equation*}
Therefore the dimension of the model space in Theorem 2.5 becomes
$\bar a= a-\epsilon^2/a$. Hence by the above argument we obtain
\begin{equation*}
\begin{split}
-\sqrt{\Lambda}r \leq  f(r)-f(0) &\leq -\sqrt{\Lambda}r +
\sqrt{n\sqrt{\Lambda}}(\sqrt{r}+1).
\end{split}
\end{equation*}

\end{document}